\documentclass[12pt, leqno]{scrartcl}

\usepackage[a4paper, hmargin=2.5cm, vmargin=3.0cm]{geometry}

\usepackage{graphicx}
\usepackage{verbatim}
\usepackage{color}
\usepackage{multicol}

\usepackage{enumerate}

\usepackage{epsf}
\usepackage[colorlinks]{hyperref}
\usepackage{booktabs}

\usepackage{tikz}
\usetikzlibrary{arrows,positioning}
\usetikzlibrary{decorations.markings}
\usetikzlibrary{arrows.meta}

\usepackage{amssymb}
\usepackage{amsmath}
\usepackage{amsfonts}
\usepackage{amsthm}
\usepackage{epstopdf}

\newcommand{\ZZ}{\mathbb{Z}}
\newcommand{\QQ}{\mathbb{Q}}
\newcommand{\qbar}{\overline{\QQ}}
\newcommand{\FF}{\mathbb{F}}

\newcommand{\CC}{\mathbb{C}}
\newcommand{\GG}{\mathbb{G}}

\newcommand{\fp}{\FF_p}
\newcommand{\fpbar}{\overline{\FF}_p}
\newcommand{\gq}{\mathrm{Gal}(\qbar/\QQ)}
\newcommand{\gm}{\GG_m}
\newcommand{\Hom}{\operatorname{Hom}}
\newcommand{\Frob}{\mathrm{Frob}}
\newcommand{\tht}{\vartheta}
\newcommand{\charpoly}{\mathrm{charpoly}}

\newcommand{\Eis}{\mathrm{Eis}}
\renewcommand{\H}{\mathrm{H}}

\newcommand{\longto}{\longrightarrow}
\newcommand{\om}{\underline{\omega}}
\newcommand{\oh}{\mathcal{O}}
\newcommand{\SL}{\operatorname{SL}}

\newcommand{\GL}{\operatorname{GL}}
\newcommand{\GSp}{\operatorname{GSp}}

\newcommand{\GU}{\operatorname{GU}}
\newcommand{\Res}{\operatorname{Res}}

\newcommand{\encircle}[1]{%
  \tikz[baseline=(X.base)]
\node (X) [draw, shape=circle, inner sep=0] {\strut #1};}

\theoremstyle{plain}
\newtheorem{theorem}{Theorem}

\theoremstyle{definition}

\theoremstyle{remark}

\newtheorem{example}[theorem]{Example}

\numberwithin{equation}{section}

\title{%
  Differential operators on modular forms (mod $p$)
}
\author{%
  Alexandru Ghitza\footnote{({\ttfamily{aghitza@alum.mit.edu}})\quad The author thanks Hidenori Katsurada for the generous financial support that made possible his participation in the Workshop on automorphic forms at RIMS in January 2018.}\\
  School of Mathematics and Statistics\\
  University of Melbourne
}
\date{\today}

\begin{document}
\thispagestyle{empty}

\maketitle
\begin{abstract}
  We give a survey of recent work on the construction of differential operators on various types of modular forms (mod $p$).
  We also discuss a framework for determining the effect of such operators on the mod $p$ Galois representations attached to Hecke eigenforms.
\end{abstract}

\section{Introduction}
\label{sect:intro}

A celebrated result of Deligne~\cite[Proposition 11.1]{Gross-tameness} attaches a mod $p$ Galois representation to a mod $p$ Hecke eigenform:
\begin{theorem}[Deligne]
  Let $f$ be a modular form (mod $p$) of weight $k$ and level $\Gamma_0(N)$ and suppose $f$ is an eigenform for all the Hecke operators: $T_\ell f = a_\ell f$ with $a_\ell\in\fpbar$.

  Then there is a continuous semisimple representation
  \begin{equation*}
    \rho_f\colon\gq\to\GL_2(\fpbar)
  \end{equation*}
  that is unramified at all primes $\ell\nmid pN$ and, for all such $\ell$, the characteristic polynomial of $\rho_f(\Frob_\ell)$ is
  \begin{equation*}
    x^2-a_\ell x+\ell^{k-1}.
  \end{equation*}
\end{theorem}

This provides the top row of a diagram that can be extended as indicated in Figure~\ref{fig:motivation}, where
\begin{itemize}
  \item the right vertical map is tensoring by the mod $p$ cyclotomic character
    \begin{equation*}
      \chi\colon\gq\to\GL_1(\fpbar),
    \end{equation*}
    determined by $\chi(\Frob_\ell)=\ell$ for all $\ell\neq p$;
  \item the effect of the left vertical map $\tht$ on $q$-expansions is simply
\begin{equation*}
  (\tht f)(q)=q\,\frac{df(q)}{dq}
\end{equation*}
The fact that the formal power series on the right hand side is once again a mod $p$ modular form (of weight $k+p+1$) is a theorem of Serre and Swinnerton-Dyer~\cite[Section 3]{SwinnertonDyer} in level one and of Katz~\cite{Katz-theta} for general level.
\end{itemize}

\begin{figure}[h]
  \begin{center}
    \begin{tikzpicture}[->,>=stealth',auto,thick]
  \node (f) at (0,3) {\parbox{4.3cm}{$f\in S_k(\Gamma_0(N);\fpbar)$\\[+10pt]$T_\ell f=a_\ell f$}};
  \node (thetaf) at (0,0) {\parbox{4.3cm}{$\tht f\in S_{k+p+1}(\Gamma_0(N);\fpbar)$\\[+10pt]$T_\ell \tht f=\ell a_\ell \tht f$}};
  \node (rhof) at (7.1,3) {\parbox{6.6cm}{$\rho_f\colon\gq\to\GL_2(\fpbar)$\\[+10pt]$\charpoly\rho_f(\Frob_\ell)=x^2-a_\ell x+\ell^{k-1}$}};
  \node (rhothetaf) at (7.1,0) {\parbox{6.6cm}{$\rho_{\tht f}\colon\gq\to\GL_2(\fpbar)$\\[+10pt]$\charpoly\rho_{\tht f}(\Frob_\ell)=x^2-\ell a_\ell x+\ell^{k+1}$}};
  \draw[|->] (f.west) to [bend right] node[label={left:$\tht$}] {} (thetaf.west);
  \draw[|->] (rhof.east) to [bend left] node[label={right:$-\otimes\chi$}] {} (rhothetaf.east);
  \draw[|->] (f) -- (rhof);
  \draw[|->] (thetaf) -- (rhothetaf);
\end{tikzpicture}
\end{center}
\caption{The relation between the theta operator on modular forms and tensoring by the cyclotomic character on Galois representations}\label{fig:motivation}
\end{figure}

The fundamental property of $\tht$ that makes the diagram possible is its commutation relation with the Hecke operators:
\begin{equation}
  \label{eq:commute_gl2}
  \tht\circ T_\ell=\ell T_\ell\circ \tht.
\end{equation}

This paper has several interrelated aims:
\begin{enumerate}[(a)]
  \item to describe various constructions of the differential operator $\tht$;
  \item to survey the generalisations of these constructions to modular forms (mod $p$) on groups other than $\GL_2/\QQ$ (in particular: Siegel, Hilbert, Picard, hermitian modular forms);
  \item to explain how commutation relations generalising Equation~\eqref{eq:commute_gl2} give rise to relations between the attached Galois representations.
\end{enumerate}

For reasons of space, we can only mention in passing the existence of similar differential operators in characteristic zero, either on $C^\infty$ (rather than holomorphic) modular forms, or on $p$-adic modular forms.
Analytic constructions of such operators go all the way back to work of Maass and of Shimura.
The modern approach uses algebraic geometric methods, which illuminate the rationality properties of these operators; for instance, see~\cite[Chapter II]{Katz-CM} for Hilbert modular forms,~\cite[Section 4]{Harris} for Siegel modular forms or~\cite{Eischen} for modular forms on unitary groups.

\subsection{Notation}
We denote the space of modular forms of weight $k$, level $\Gamma_0(N)$ and coefficient ring $R$ by $M_k(N;R)$.
The cuspidal subspace will be denoted $S_k(N;R)$.

\section{Analytic construction}
\label{sect:analytic}
We begin by describing this construction in the case of classical modular forms.
Up to minor modifications, it is a special case of the method employed by B\"ocherer and Nagaoka in the Siegel setting~\cite[Section 4]{BoechererNagaoka-firsttheta}.

It was observed long ago that the derivative of $f\in M_k(N;\CC)$ satisfies the not-quite-modular relation
\begin{equation*}
  f^\prime\left(\frac{az+b}{cz+d}\right)
  =(cz+d)^{k+2}f^\prime(z)+ck(cz+d)^{k+1}f(z).
\end{equation*}
This failure to be modular can be exploited/fixed in a number of ways, for which we refer the reader to Zagier's expositions in~\cite[Chapter 5]{Zagier-123},~\cite{Zagier-diff}.

One of these approaches involves taking a second form $g\in M_\ell(N;\CC)$ and defining the Rankin-Cohen bracket
\begin{equation*}
  [f,g]=\frac{1}{2\pi i}\left(kfg^\prime-\ell f^\prime g\right)\in M_{k+\ell+2}(N;\CC).
\end{equation*}

This construction preserves the ring of coefficients of the Fourier expansions: if $R$ is a subring of $\CC$ and both $f$ and $g$ have coefficients in $R$, then so does $[f,g]$.

We can use the Rankin-Cohen bracket to construct a theta operator on mod $p$ classical modular forms, as described in the following diagram:
\begin{center}
  \begin{tikzpicture}[->,>=stealth',auto,thick]
    \node (fbar) at (0,0) {$\overline{f}\in S_k(N;\fpbar)$};
    \node (f) at (0,3) {$f\in S_k(N;\ZZ[1/N])$};
    \node (bracket) at (7,3) {$[f,g]\in S_{k+p+1}(N;\ZZ[1/N])$};
    \node (thetaf) at (7,0) {$\tht(\overline{f}):=\overline{[f,g]}\in S_{k+p+1}(N;\fpbar)$};
    \draw[|->] (fbar) -- node [label=left:{\encircle{L}},label=right:{if $k\geq 2$}] {} (f);
    \draw[|->] (f) -- node [label=above:{\encircle{B}}] {} (bracket);
    \draw[|->] (bracket) -- node [label=left:{mod $p$},label=right:{\encircle{R}}] {} (thetaf);
    \draw[|->,dotted] (fbar) -- (thetaf);
\end{tikzpicture}
\end{center}

In step~\encircle{L} we lift the mod $p$ cusp form $\overline{f}$ to a form $f$ in characteristic zero.
This is known to be possible if the weight is at least two~\cite[Theorem~1.7.1]{Katz-modular}.
In step~\encircle{B} we take the Rankin-Cohen bracket of $f$ and an auxiliary form $g$ such that $g(q)\equiv 1\pmod{p}$; concretely, we use the Eisenstein series $g=E_{p-1}$ provided that $p\geq 5$.
Finally, step~\encircle{R} is simply reduction modulo $p$.

The ingredients involved in this construction have been generalised to other settings:
\begin{enumerate}[(a)]
  \item Lifting modular forms mod $p$ to characteristic zero is possible in many situations if both the weight and the prime $p$ are sufficiently large, as shown for PEL-type Shimura varieties by Lan-Suh~\cite{LanSuh-lifts},~\cite[Theorem 8.13 and Corollary 8.14]{LanSuh-vanishing}, and in the special case of Siegel modular forms by Hida~\cite[Section 3.5]{Hida} and Stroh~\cite{Stroh-siegel},~\cite[Th\'eor\`eme 1.3]{Stroh-classicite}.
    In cases not covered by these results, one can dispense with the lifting step at the expense of getting an operator that is only defined on the subspace of modular forms mod $p$ that are liftable\footnote{These are colloquially referred to as mod $p$ forms \emph{\`a{} la Serre} (defined as reductions modulo $p$ of characteristic zero forms), as opposed to mod $p$ forms \emph{\`a{} la Katz} (defined intrinsically as sections of line bundles on modular curves in positive characteristic).}.
  \item Rankin-Cohen brackets have been defined for Siegel modular forms~\cite[Theorem 1.4]{ChoieEholzer},~\cite{EholzerIbukiyama}, hermitian modular forms~\cite[Theorem 3.2]{MartinSenadheera}, Hilbert modular forms~\cite[Theorem 5.1]{Lee},~\cite[Section 4]{ChoieKimRichter}, Picard modular forms~\cite[Section 9]{CleryVanDerGeer},~\cite{FreitagSalvatiManni}.
    There are typically several choices of brackets, many of which start with scalar-valued forms and yield vector-valued forms\footnote{Finding generators for spaces of vector-valued forms is in fact a major application of Rankin-Cohen brackets.}.
  \item Characteristic zero forms whose $q$-expansion is congruent to $1$ modulo $p$: the state of the art is the construction of such forms as theta series attached to certain special lattices, which was accomplished in the case of Siegel modular forms~\cite[Theorem 1]{BoechererNagaoka-firsttheta} and hermitian modular forms~\cite[Theorem 3.3]{HentschelNebe},~\cite[Theorem 3.1]{KikutaNagaoka-congruence},~\cite[Propositions 5.1 and 6.1]{KikutaNagaoka-modp}.
    They are obviously related to the Hasse invariants that appear in the algebraic-geometric context described in the next section.
\end{enumerate}

Putting everything together, there are analytic constructions of theta operators for Siegel modular forms~\cite[Theorem 4]{BoechererNagaoka-firsttheta},~\cite[Theorem 4.5]{BoechererNagaoka-secondtheta} and hermitian modular forms~\cite[Theorem 3]{KikutaNagaoka-theta}.

A drawback of the method described in this section is that it is by no means clear how the resulting theta operators commute with the Hecke operators.
This needs to be checked separately by means of rather unenlightening calculations involving explicit formulas for the effect of Hecke operators on $q$-expansions (see for example~\cite[Corollary 15]{theta-galois} for the case of the B\"ocherer-Nagaoka theta operators on Siegel modular forms).

\section{Algebraic-geometric construction}
\label{sect:geometric}
This originated with Katz~\cite{Katz-theta}; we will describe a variant due to Gross~\cite[Section 5]{Gross-tameness}, that was first brought to our attention by Eyal Goren.

Let $X^{\text{ord}}$ denote the ordinary locus of the modular curve $X:=X_0(N)_{\fpbar}$ and consider the Igusa curve $\tau\colon I\to X^{\text{ord}}$.

The modular curve $X$ is endowed with an invertible sheaf $\om$ (the Hodge bundle) which gives rise to the spaces of modular forms:
\begin{equation*}
  M_k(N;\fpbar)=\H^0\left(X,\om^{\otimes k}\right).
\end{equation*}

The construction proceeds as indicated in Figure~\ref{fig:geometric}.

\begin{figure}[h]
\begin{center}
  \begin{tikzpicture}[->,>=stealth',auto,thick]
    \node (f) at (0,6) {$f\in\H^0(X,\om^{\otimes k})$};
    \node (tauf) at (0,4) {$\tau^* f\in\H^0(I,\tau^*\om^{\otimes k})$};
    \node (taufak) at (0,2) {$\frac{\tau^* f}{a^k}\in\H^0(I,\oh_I)$};
    \node (etaf) at (0,0) {$\eta_f:=d\left(\frac{\tau^* f}{a^k}\right)\in\H^0(I,\Omega_I^1)$};
    \node (kappa) at (10,0) {$\kappa^{-1}(\eta_f)\in\H^0(I,\tau^*(\om^{\otimes 2}(-C)))$};
    \node (akkappa) at (10,2) {$a^k\kappa^{-1}(\eta_f)\in\H^0(I,\tau^*(\om^{\otimes k+2}(-C)))$};
    \node (g) at (10,4) {$g\in\H^0(X^{\text{ord}},\om^{\otimes k+2}(-C))$};
    \node (thetaf) at (10,6) {$\tht f:=Ag\in\H^0(X,\om^{k+p+1}(-C))$};
    \draw[|->] (f) -- (tauf);
    \draw[|->] (tauf) -- (taufak);
    \draw[|->] (taufak) -- node [label=left:{$d\colon\oh_I\to\Omega_I^1$}] {} (etaf);
    \draw[|->] (etaf) -- node [label=below:{$\kappa\colon\tau^*(\om^{\otimes 2}(-C))\to\Omega_I^1$}] {} (kappa);
    \draw[|->] (kappa) -- (akkappa);
    \draw[|->] (akkappa) -- (g);
    \draw[|->] (g) -- node [label=right:{\encircle{E}}] {} (thetaf);
    \draw[|->,dotted] (f) -- (thetaf);
\end{tikzpicture}
\end{center}
\caption{Algebraic-geometric construction of theta operator on classical modular forms mod $p$: the Igusa curve version}
\label{fig:geometric}
\end{figure}

The lower horizontal map is induced by the Kodaira-Spencer isomorphism
\begin{equation*}
\kappa\colon\tau^*(\om^{\otimes 2}(-C))\to\Omega_I^1,
\end{equation*}
where $C$ is the divisor consisting of the cusps on the modular curve $X$.

Finally, in step~\encircle{E} we extend the section $g$ from the ordinary locus to the entire modular curve: a calculation shows that $g$ has a simple pole at each supersingular point in $X$, so multiplying $g$ by the Hasse invariant $A\in\H^0(X,\om^{\otimes p-1})$ clears the poles.

The commutation relation~\eqref{eq:commute_gl2} between $\tht$ and the Hecke operator $T_\ell$ can be obtained by determining the effect of degree $\ell$ isogenies on the various steps.
All steps commute with $\ell$-isogenies, with the exception of the Kodaira-Spencer isomorphism, which commutes up to multiplication by $\ell$.

The original form of this argument, given in~\cite{Katz-theta}, works entirely on the modular curve $X$.
The role of differentiation is played by the Gauss-Manin connection $\nabla$, which is defined on de Rham cohomology (instead of sheaf cohomology).
Hence it is necessary to pass from sheaf cohomology to de Rham cohomology and back, which can be done over the ordinary locus of $X$.
Extending to all of $X$ once again involves multiplication by the Hasse invariant.
For the details, we refer the interested reader to~\cite{Katz-theta} or the exposition in~\cite{Ramsey-theta}.

Some of the ingredients in the above constructions are automatically defined in much greater generality.
Others have been extended to more general settings:
\begin{enumerate}[(a)]
  \item Igusa varieties for certain unitary groups~\cite[Chapter IV]{HarrisTaylor} or for more general PEL-type Shimura varieties~\cite[Section 4]{Mantovan-pel};
  \item Hasse invariants: Hilbert case covered in~\cite[Section 7]{AndreattaGoren}, see~\cite{KoskivirtaWedhorn,GoldringKoskivirta} for recent developments for general Shimura varieties.
\end{enumerate}

Putting these pieces together\footnote{We took an intentionally over-simplistic view of the process, in order to allow the big picture of the algebraic-geometric framework to emerge.
In practice, each setting has its own geometric intricacies, and  dealing with these requires more than just piecing together the various ingredients.}, there are algebraic-geometric constructions of mod $p$ theta operators for
\begin{enumerate}[(a)]
  \item Hilbert modular forms~\cite[Section 15]{AndreattaGoren};
  \item Siegel modular forms in genus two~\cite[Proposition 3.9]{Yamauchi} and general genus~\cite{FlanderGhitza};
  \item Picard modular forms~\cite[Section 3]{deShalitGoren-picard}, more generally forms on $\GU(m,n)$~\cite{deShalitGoren-unitary} and~\cite[Sections 5 and 6]{EischenMantovan}.
\end{enumerate}

\section{Group-cohomological construction}
\label{sect:cohomological}
For a ring $R$, let $V_{k-2}(R)=R[x,y]_{k-2}$ denote the space of homogeneous polynomials of degree $k-2$.
We consider the cohomology group $\H^1(\Gamma,V_{k-2}(R))$ where $\Gamma$ is a congruence subgroup of $\SL_2(\ZZ)$.

The link with modular forms comes from the Eichler-Shimura isomorphism
\begin{equation*}
  \H^1(\Gamma, V_{k-2}(\CC))\cong S_k(\Gamma;\CC)\oplus
  \overline{S_k(\Gamma;\CC)}\oplus\Eis_k(\Gamma;\CC),
\end{equation*}
where the bar indicates complex conjugation and $\Eis_k$ is the space spanned by the Eisenstein series of weight $k$.

We now focus on the case where the ring of coefficients $R=\fp$.
Consider the element $\theta(x,y)=x^py-xy^p\in V_{p+1}(\fp)$.
Multiplication by $\theta$ induces a map in cohomology
\begin{equation*}
  \H^1(\Gamma,V_{k-2}(\fp))\to\H^1(\Gamma,V_{k+p+1-2}(\fp))
\end{equation*}
that commutes with the Hecke operator $T_\ell$ up to multiplication by $\ell$.
The point is that $\theta$ is a polynomial semi-invariant of the finite group $\GL_2(\fp)$.

Despite the fact that the relation between group cohomology and modular forms is less clear-cut in the case of mod $p$ coefficients than for complex coefficients, there is a Hecke action and there is value in studying Hecke eigenclasses in group cohomology in their own right.
The foundations were laid by Ash and Stevens in~\cite{AshStevens-cohomology} and pursued by them in the $\GL_2$ case in~\cite{AshStevens-modular}, which includes the description of the group cohomological theta operator we gave above.

There has been a significant amount of work done in generalising this to higher $\GL_n$, mostly by Ash and his collaborators\footnote{At the risk of omitting names, these include Doud, Gunnells, McConnell, D.~Pollack, Sinnott, Stevens.}.
This was motivated by their success in performing explicit computations\footnote{For an exposition of the computational methods used to study the cohomology of arithmetic groups such as $\GL_n$, see~\cite{Gunnells}.} with $\GL_3$ and $\GL_4$, but also by the fact that, in sharp contrast to the classical case $\GL_2$, for $n>2$ there is no direct link between automorphic representations and the algebraic-geometric framework of Shimura varieties that plays such a central role in the Langlands program.
Indeed, they formulated precise conjectures positing the existence of $n$-dimensional mod $p$ Galois representations attached to Hecke eigenclasses in the cohomology of $\GL_n$ with $\fp$-coefficients, see~\cite{Ash-glnz,Ash-modp,AshDoudPollack}.
Recent breakthroughs due to Scholze resulted in the proof of a large part of these conjectures; we refer the interested reader to~\cite{Scholze-torsion}, as well as the survey~\cite{Morel} and the recent improvements~\cite{Caraiani,NewtonThorne}.

There are group-cohomological approaches to the study of modular forms on other groups as well, for instance:
\begin{itemize}
  \item Hilbert modular forms~\cite{Reduzzi-shiftings} (dedicated to ``weight shiftings'', the cohomological analogues of theta operators), see also the computational framework in~\cite{GunnellsYasaki};
  \item Bianchi modular forms, see~\cite{SengunTurkelli,Sengun-integral,Sengun,Gunnells-bianchi};
  \item Siegel modular forms~\cite{Buecker}.
\end{itemize}

\section{Modular forms on other groups}
\label{sect:modular}
We return to the task of generalising Figure~\ref{fig:motivation} to groups other than $\GL_2/\QQ$.
In order to describe the effect of (some of) the operators defined in the previous sections on the corresponding Galois representations, we need to introduce some of the representation-theoretic foundations of modular forms on algebraic groups.
To keep the exposition clean, we restrict our attention to groups over $\QQ$; the correct level of generality is to work over number fields $F$, and we invite interested readers to consult~\cite{TreumannVenkatesh,BuzzardGee} for the more general setup.

Let $G$ be a connected reductive algebraic group over $\QQ$.
We will take it for granted that there is a notion of modular form (mod $p$) on the group $G$.
This may, for instance, come from a Shimura variety attached to $G$, and is known to be the case for the various specific types of modular forms we consider in this paper, as summarised in Table~\ref{table:forms}.

Fixing subgroups $T\subset B\subset G$ where $T$ is a maximal torus and $B$ is a Borel subgroup, we have the root datum
\begin{equation*}
  \left(X^\bullet, \Delta^\bullet, X_\bullet, \Delta_\bullet\right),
\end{equation*}
where
\begin{itemize}
  \item $X^\bullet=\Hom(T,\gm)$ is the group of characters of $T$;
  \item $X_\bullet=\Hom(\gm,T)$ is the group of cocharacters of $T$;
  \item $\Delta^\bullet$ is the set of simple roots, with corresponding simple coroots $\Delta_\bullet$.
\end{itemize}

There is a natural pairing
\begin{equation*}
  \langle\cdot,\cdot\rangle\colon X^\bullet\times X_\bullet\longto \Hom(\gm,\gm)=\ZZ
\end{equation*}
given by composition: $\langle \alpha,\lambda\rangle=k$ where $(\alpha\circ\lambda)(z)=z^k$ for all $z\in\gm$.

The dual group $\hat{G}$ has subgroups $\hat{T}\subset\hat{B}\subset\hat{G}$ and root datum
\begin{equation*}
  \left(X_\bullet, \Delta_\bullet, X^\bullet, \Delta^\bullet\right).
\end{equation*}
In other words, the character $\alpha\colon T\to\gm$ of $T$ is identified with the cocharacter $\hat{\alpha}\colon\gm\to \hat{T}$ of $\hat{T}$.

\begin{table}[h]
  \begin{center}
  \begin{tabular}{ll}
    \toprule
    type of modular form & algebraic group\\ 
    \midrule
    classical & $\GL_2$\\ 
    Siegel & $\GSp_{2g}$\\ 
    Hilbert & $\Res_{F/\QQ}\GL_2$, $F$ totally real\\ 
    Bianchi & $\Res_{K/\QQ}\GL_2$, $K$ imaginary quadratic\\ 
    Picard & $\GU_{2,1}$\\ 
    hermitian & $\GU_{n,n}$\\ 
    \bottomrule
  \end{tabular}
  \end{center}
  \caption{Some types of modular forms with the corresponding algebraic groups over $\QQ$}
  \label{table:forms}
\end{table}

Given a prime $\ell$, the local Hecke algebra at $\ell$ is
\begin{align*}
  \mathcal{H}_\ell(G)=\mathcal{H}(G(\QQ_\ell), G(\ZZ_\ell))
  =\{&t\colon G(\QQ_\ell)\to\ZZ\text{ locally constant, compactly}\\
     &\text{supported and $G(\ZZ_\ell)$-bi-invariant}\}.
\end{align*}
This is a free $\ZZ$-module with a basis of characteristic functions
\begin{equation*}
  t_{\ell,\lambda} = \mathrm{char}\left(G(\ZZ_\ell)\lambda(\ell)G(\ZZ_\ell)\right)
\end{equation*}
indexed by $\lambda\in P^+$, the cone of dominant coweights of $G$:
\begin{equation*}
  P^+=\left\{\lambda\in X_\bullet\mid \langle\alpha,\lambda\rangle \geq 0\text{ for all }\alpha\in\Delta^\bullet\right\}.
\end{equation*}
Since $P^+$ is identified with the dominant weights of $\hat{G}$, it also indexes the irreducible representations of $\hat{G}$.

This relation is made more explicit by the Satake isomorphism~\cite[Proposition~3.6]{Gross-satake}
\begin{equation*}
  \mathcal{S}_{\ZZ,\ell}\colon \mathcal{H}_\ell(G)\otimes\ZZ[\ell^{\pm 1/2}]\longto R(\hat{G})\otimes\ZZ[\ell^{\pm 1/2}].
\end{equation*}
We tensor this with $\fpbar$ for $p\neq\ell$ to get an isomorphism
\begin{equation*}
  \mathcal{S}_\ell\colon\mathcal{H}_\ell(G)\otimes\fpbar\longto R(\hat{G})\otimes\fpbar.
\end{equation*}

Suppose that we are given a Hecke module, that is a finite-dimensional $\fpbar$-vector space $V$ endowed with an action of $\prod_{\ell\neq p}\mathcal{H}_\ell(G)$.
(The examples we have in mind are spaces of modular forms (mod $p$) attached to the group $G$.)
Given a Hecke eigenform $f\in V$, we get the commutative diagram
\begin{center}
  \begin{tikzpicture}[->,>=stealth',auto,thick]
    \node (h) at (0,2) {$\mathcal{H}_\ell(G)\otimes\fpbar$};
    \node (f) at (0,0) {$\fpbar$};
    \node (r) at (4,2) {$R(\hat{G})\otimes\fpbar$};
    \draw[->] (h) -- node [label={$\mathcal{S}_\ell$}] {} (r);
    \draw[->] (h) -- node [label={left:$\Psi_{f,\ell}$}] {} (f);
    \draw[->,dotted] (r) -- node [label={right:$\omega_f$}] {} (f);
\end{tikzpicture}
\end{center}
where
\begin{itemize}
  \item the vertical arrow $\Psi_{f,\ell}\colon\mathcal{H}_\ell(G)\otimes\fpbar\to\fpbar$ is the local Hecke eigensystem attached to $f$, that is the algebra homomorphism that maps a Hecke operator $t$ to its eigenvalue: $tf=\Psi_{f,\ell}(t)f$ for all $t\in\mathcal{H}_\ell(G)\otimes\fpbar$;
  \item the diagonal arrow $\omega_f\colon R(\hat{G})\otimes\fpbar\to\fpbar$ is the character of the representation group defined by $\omega_f:=\Psi_{f,\ell}\circ\mathcal{S}_\ell^{-1}$.
\end{itemize}

But characters of $R(\hat{G})\otimes\fpbar$ are in bijective correspondence with semisimple conjugacy classes $s$ in $\hat{G}(\fpbar)$, via $\omega_s(\lambda)=\lambda(s)$ for all $\lambda\in X_\bullet$.
So the character $\omega_f$ determines a semisimple conjugacy class $s_{f,\ell}\in\hat{G}(\fpbar)$, called the Satake parameter of $f$ at $\ell$.

Under the mod $p$ Langlands correspondence, there is a general expectation that, to a mod $p$ Hecke eigenform $f$ on the group $G$, one can associate a mod $p$ Galois representation $\rho_f\colon\gq\to\hat{G}(\fpbar)$ such that $\rho_f(\Frob_\ell)=s_{f,\ell}$ for all but finitely many primes $\ell\neq p$.

This gives the natural generalisation of the top row of the diagram in Figure~\ref{fig:motivation}; the next section aims to generalise the rest of the diagram.

\section{Effect on Galois representations}
\label{sect:galois}

Given a Hecke eigenform $f$ defined on some group, and its image $\tht f$ under a theta operator of the type discussed in the previous sections, how does one relate the Galois representation attached to $\tht f$ to that attached to $f$?
In the classical case illustrated in Figure~\ref{fig:motivation}, this follows easily from comparing the characteristic polynomials of Frobenius elements.
In general, such a direct calculation is not feasible, but the relation between Galois representations can still be deduced, as an almost formal consequence of the commutation relation between $\tht$ and the Hecke operators.

More precisely, using the Satake isomorphism we can prove (with the notation from Section~\ref{sect:modular}):
\begin{theorem}[\cite{theta-galois}]
  \label{thm:galois}
  Let $f$ be a Hecke eigenform (mod $p$) on a group $G$.
  Let $\tht$ be a map of modular forms such that $\tht f$ is a Hecke eigenform whose Hecke eigensystems satisfy
  \begin{equation*}
    \Psi_{\tht f,\ell}(t_{\ell,\lambda})=\eta(\lambda(\ell)) \Psi_{f,\ell}(t_{\ell,\lambda})\qquad
    \text{for all }\ell\notin\Sigma,\lambda\in P^+
  \end{equation*}
  where $\eta\colon G\to\gm$ is a character of $G$ and $\Sigma$ is a finite set of primes containing $p$.
  Then
  \begin{equation*}
    \rho_{\tht f}\cong (\hat{\eta}\circ\chi) \rho_f
  \end{equation*}
  as Galois representations $\gq\to\hat{G}(\fpbar)$, where $\chi$ is the mod $p$ cyclotomic character.
\end{theorem}

We illustrate this in the simplest setting:

\begin{example}[Classical modular forms]
  Take $G=\GL_2/\QQ$.

  There are two dominant coweights $\lambda_i\colon\gm\to T$ to consider:
  \begin{equation*}
    \lambda_1\colon z\mapsto \begin{bmatrix}1&0\\0&z\end{bmatrix},
    \qquad
    \qquad
    \lambda_2\colon z\mapsto \begin{bmatrix}z&0\\0&z\end{bmatrix}.
  \end{equation*}
  The corresponding elements of the local Hecke algebra are the Hecke operator $T_\ell$ and multiplication by $\ell$.

  We take $\eta:=\det\colon G\to\gm$, and we verify easily that (over $\fpbar$) we have
  \begin{equation*}
    a_\ell(\tht f) = \ell a_\ell(f),
    \qquad
    \qquad
    \ell^{(k+p+1)-1}=\ell^2 \ell^{k-1}.
  \end{equation*}
  
  Since $\hat{\eta}\colon\gm\to\hat{G}\cong G$ is given by
  \begin{equation*}
    \hat{\eta}\colon z\mapsto \begin{bmatrix}z&0\\0&z\end{bmatrix},
  \end{equation*}
  Theorem~\ref{thm:galois} says that
  \begin{equation*}
    \rho_{\tht f}\cong \begin{bmatrix}\chi&0\\0&\chi\end{bmatrix}\,\rho_f\cong\chi\otimes\rho_f,
  \end{equation*}
  recovering what we already knew from the direct comparison of the characteristic polynomials of $\Frob_\ell$.
\end{example}

\bibliographystyle{alpha}
\bibliography{refs}

\end{document}